\documentclass[10pt]{amsart}
\usepackage{rotating}
\usepackage{amssymb}
\usepackage{url}
\usepackage[all]{xy}
\usepackage{color}
\usepackage[table]{xcolor}
\definecolor{orange}{rgb}{1,.7,0}
\definecolor{violet}{rgb}{.5,0,.5}
\definecolor{dg}{rgb}{0,0.67,0}
\definecolor{g}{rgb}{.75,.75,.75}
\definecolor{lg}{rgb}{.90,.90,.90}
\newcommand{\Z}{{\mathbb Z}}

\newcommand{\Fr}{{\mbox{Fr}}}

\newcommand{\Q}{{\mathbb Q}}
\newcommand{\F}{{\mathbb F}}
\newcommand{\C}{{\mathbb C}}

\newcommand{\Gal}{\mbox{Gal}}
\newcommand{\ord}{\mbox{ord}}

\newcommand{\cmmt}[1]{}
\newtheorem{Theorem}{Theorem}[section]

\allowdisplaybreaks

\title[$PGL_2(\F_\ell)$ number fields with rational companion forms]{$PGL_2(\F_\ell)$ number fields \\ with rational companion forms}

\author{David P.\ Roberts}
\address{Division of Science and Mathematics, University of
  Minnesota-Morris, Morris, MN 56267}
\email{roberts@morris.umn.edu}
\begin{document}
\maketitle 
\begin{abstract}
We give a list of $PGL_2(\F_\ell)$ number fields for $\ell \geq 11$ 
which have rational companion forms.
Our list has fifty-three fields and seems likely to be complete.   
Some of the fields on our list are very lightly ramified for their 
Galois group.  
%
\end{abstract}


\section{Introduction}
   One of the great recent advances in number theory is the proof by
 Khare and Wintenberger \cite{kw1,kw2}  of the Serre reciprocity conjecture \cite{se}.   The Khare-Wintenberger reciprocity theorem says that 
 certain two-dimensional Galois representations in prime characteristic $\ell$ necessarily come from classical modular forms.  
 Once certain normalizations are in place, a given representation $\rho$ typically
 comes from just one form $g$.   
  However if $\rho$ is tamely or minimally wildly 
  ramified at $\ell$,
  then it comes from 
 two forms $g$ and $h$, called companions.  
 
 The purpose of this paper is to give
 a systematic collection of fifty-three examples of this phenomenon of companion
 forms.     Our collection illustrates all three types of companion forms, which we 
  label 1T or {\em diagonalizable},
 2T or {\em supersingular}, and  2W or {\em peu ramifi\'ee}. 
 Important early theoretical developments on companion
 forms in these senses took place in letters among 
 Deligne, Fontaine, and Serre in the 1970s.  The theory for Type 1T was established by Gross in 1990 \cite{gr}.  
 Edixhoven gave the first complete proofs for cases 2T and 2W in
 1992 \cite{Edi92}.   To be noted is that the term ``companion form'' is often 
 used in the context of 1T only.  We are using it more broadly, because
 the three cases are very similar from the viewpoint of this paper. 


 We aim to address as general a readership as possible.  Thus, to the extent possible,
 we work with explicit degree $\ell+1$ polynomials  with Galois group $PGL_2(\F_\ell)$ 
 rather than Galois representations.   Similarly, it entirely suffices for us to regard modular forms as 
 power series in a formal variable $q$.  

 
 In Section~\ref{two}, we present two of the fifty-three examples in some detail,
 both of diagonalizable type with $\ell=11$.  
 Each example starts from a degree twelve polynomial $f(x) \in \Q[x]$
 with Galois group $PGL_2(\F_{11})$ and ends with a pair $(g,h)$ 
 of companion forms in $\Q[[q]]$.    Our focus is not so much on the 
 polynomials themselves, but rather on the degree twelve fields
 $K = \Q[x]/f(x)$ they define.  
 
  In Section~\ref{fifty-three}, we ask for all fields $K$ belonging to triples $(K,g,h)$ of the same general nature, 
 including the very strong requirement that both $g$ and $h$ have
 coefficients in $\Q$. We restrict to $\ell \geq 11$ to keep the final collection
 of manageable size.  Pairwise comparing all known
 rational forms, we extract those pairs $(g,h)$ which
 satisfy the companionship condition.  We find
 in Theorem~\ref{mt} fifty-three fields $K$ which are part of 
 triples in this way, and we conjecture that there are no more.    
 The distribution of the number fields $K$ with regard to $\ell$ and the 
 three 
 notions of companionship is as follows: 
  \[
  \begin{array}{r    |cccccccccc}
       \;\;\;\;\;\;\;\;\;\;\;\;\;\;\;\;\;\; \ell: &  11 & 13 & 17 & 19 & 23 & 29 & 31 & 37 & 41 & 43  \\
   \hline
  \mbox{1T=Diagonalizable:} & 6 & 5 & 2 & 2 & 1 & 2 & 1 &  & 1 &   \\
 \mbox{2T = Supersingular:}  & 8 & 7 &4 & 3 &  & 3 &  & 1 &   &  1 \\
 \mbox{2W=Peu ramifi\'ee:} & 5 &  & 1 &  &  &  &  \\
  \end{array}.
  \]
 While
 the forms themselves can easily be made very explicit,
our computations in Section~\ref{fifty-three} do not produce defining polynomials $f(x)$ for the fields $K$.

  In Section~\ref{lightly}, we study the fifty-three fields further, finding defining 
  polynomials when possible.  Some of these number fields $K$ 
  are very lightly ramified for their Galois group.   
   We explain how incorporating our fields into a systematic tabulation of number 
  fields based on Serre reciprocity seems possible.  This systematic tabulation would 
  involve totally dropping our rationality conditions.   To a large extent, it 
  would then consist of repeating our computations here in the resulting larger context.   
  In particular, beating our best fields or establishing them as true minima
  seems within reach by modular methods.
  
  Computations in this paper were done using a mix of {\em Magma} \cite{Magma},
  {\em Pari} \cite{gp}, and {\em Mathematica} \cite{mma}.   Together with
  the closely related paper \cite{RobNRC}, this paper grew
  from a talk given by the author at {\em Automorphic Forms: theory and computation} at King's College
  London, in September 2016.  The author's research was supported by grant \#209472 from
  the Simons Foundation and grant DMS-1601350 from the National Science Foundation.    
  

 \section{Two $PGL_2(\F_{11})$ fields with rational companion forms}   
 \label{two}

    In this section we present two remarkably parallel examples,
centering on triples $(K_1,g_1,h_1)$ and $(K_2,g_2,h_2)$.  
We exhibit very concretely how Galois-theoretic invariants of the $K_i$ 
are connected with modular invariants from the $(g_i,h_i)$.  
We keep background theory to an absolute minimum,
with some of this theory being presented in the next two sections.
 \subsection{Mathieu group sources}  
 \label{mathieu}
  Our two examples
  have the added interest that the number fields were first found ``accidentally''
 in a context very far removed from elliptic curves and modular forms.
  Let
 \begin{eqnarray}
\label{examp1} \;\;\;\; f_1(x) & = &   x^{12}-4 x^{11}-4 x^{10}+16 x^9+24 x^8-30 x^7 \\ 
\nonumber && \qquad  -78 x^6-18 x^5+72 x^4+86 x^3+52 x^2+16 x+2, \\
\nonumber&&\\
\label{examp2} \;\;\;\; f_2(x) & = &   x^{12}-6 x^{10}-6 x^9-6 x^8+126 x^7+104 x^6 \\ 
\nonumber && \qquad -468 x^5 +258 x^4+456 x^3-1062 x^2+774 x-380.
 \end{eqnarray}
 The fact that both polynomials have Galois group $PGL_2(\F_{11})$ 
 can be rapidly confirmed by {\em Magma}'s \verb@GaloisGroup@. 
 Their exotic source in each case is related to one of the five sporadic
 simple groups $M_n$ discovered by Mathieu in the mid-1800s.  The first comes
 from the degenerate specialization at $y = -47^2/2^6 3$ of Malle's one-parameter 
 family of $M_{22}.2$ fields \cite[Theorem~2]{ma}.  The second comes from 
 the degenerate specialization at $t=-17^3/2^7$ of a one-parameter family 
 of $M_{12}.2$ fields \cite[Cover~D2, Table~4.5]{ro-mathieu}.  These two sources actually 
 give polynomials of degree $22$ and $24$ respectively, and
 we used {\em Magma}'s \verb@GaloisSubgroup@ to 
 obtain degree twelve polynomials with the same splitting
 field.   Finally we used {\em Pari}'s \verb@polredabs@ to 
 reduce the size of coefficients.  
 
 
 \subsection{Frobenius partitions}   
 \label{frobenius}
 For a degree $n$ number field $K = \Q[x]/f(x)$ 
 and a prime $p$ not dividing its discriminant $D$, one has 
 a factorization partition $\lambda_p$.  The parts
 of $\lambda_p$ are the degrees of the irreducible
 factors of $f(x)$ in $\Q_p[x]$.  If $p$ does not divide
 the discriminant $Dc^2$ of the defining polynomial $f(x) \in \Z[x]$,
 then $\lambda_p$ is more easily computed as the 
 degrees of the irreducible factors of $f(x)$ 
 in $\F_p[x]$.     
 
 The polynomial discriminants of $f_1$ and $f_2$ are respectively
 \begin{align*}
 D_1 c_1^2 & = - 2^{14} \cdot 3^{30} \cdot 11^9, &
 D_2 c_2^2 & = -2^{12} \cdot  3^{14} \cdot 11^9  \cdot 17^2 \cdot  1907473^2 \cdot 2615189^2.
 \end{align*}
 From the construction of the stem fields $K_i=\Q[x]/f_i(x)$ via specialization
 of understood covers, one knows 
that in each case $2$, $3$, and $11$ are the primes dividing
the field discriminant $D_i$.   Table~\ref{frobtable} gives
 in each case the factorization partition $\lambda_p$ for the twenty-two 
good primes $p$ less than a hundred.   It also gives the parity 
$d_p$ of the partition.  Both polynomial discriminants are $-11$ modulo 
squares, so the $d_p$ agree in the two cases.   More explicitly,
$p$ occurs on a row with $d_p=+$ 
if and only if $p$ is a square modulo $11$, i.e.\ if and only
if $p \equiv 1, 3, 4, 5, 9 \; (11)$.  

 \begin{table}[htb]
\[
\begin{array}{|ccc|c|   l  |  l  |c|}
\hline
                    & \lambda_p &  d_p & \mbox{mass} &  \mbox{Primes $p$ for $K_1$} & \mbox{Primes $p$ for $K_2$} &
                    s_p \\
\hline
\mbox{inert}&12 & - & 1/6 & 13, 29, 79, 83 & 7, 13, 61,73, 83 & 7,8 \\
\mbox{torus}&4^3 & - & 1/12 & 7,17 & 29 &  2  \\
&6^2 & + & 1/12 &  31, 47, 97 & 23,59 & 3 \\
&3^4 & + &1/12 & 37, 71 & 37, 47 &  1   \\
&2^6 & + &1/24 &  & 71 &  0 \\
\hline
\mbox{split}&10^1 1^2 & - & 1/5 & 19,41,43,61,73  & 19,41,43,79 & 6,10 \\
\mbox{torus}&2^51^2 & - & 1/20 &  & 17 & 0  \\
       &5^21^2 & + & 1/5 &  5,23,53,59,67,89 & 53,67,89,97 & 5,9 \\ 
\hline
\mbox{uni-}&11^1 1 & + & 1/11 &   & 5,31 &  4 \\
\mbox{potent} &1^{12} & + & 1/1320 &   &   &   4 \\
\hline
\end{array}
\]
\caption{\label{frobtable} Frobenius partitions $\lambda_p$ correlating with modular quantities 
$s_p = a_p^2/p^{k-1}$ 
in our two examples. 
}
\end{table}

In general, a Frobenius partition $\lambda_p$ reflects a more refined invariant,
a conjugacy class $\Fr_p$ in the Galois group.  For the 
group $PGL_2(\F_{11}) \subset S_{12}$, there are $13$ conjugacy classes.  They
give rise to ten of the
seventy-seven partitions of twelve, and these ten partitions are given in 
Table~\ref{frobtable}.  For context, the column ``mass'' gives the
 asymptotic frequency of each partition.
 Frobenius classes $\Fr_p$ and 
 partitions behave similarly for all $PGL_2(\F_\ell)$, and the words
 in the first column summarize a structure theory that applies for all $\ell$.


 \subsection{Ramification} 
 \label{ramification}
 The discriminants of the two fields
 $K_i$ can be calculated directly, say by {\em Pari}'s \verb@nfdisc@:
 \begin{align}
 D_1 & = -2^{14} \cdot 3^{10} \cdot 11^9, &
 D_2 & = -2^{10} \cdot 3^{14} \cdot 11^9. 
 \end{align}
 It is important for us to have a clear picture of the inertia
 groups $I_p \subset PGL_2(\F_\ell)$  underlying the discriminants.  This information is given automatically
 by the $p$-adic identifier at the website of \cite{jr-local-database}.  Here, underlying
 the two exponents $10$, ramification is tame of order $|I_p|=11$.  
 Underlying the two instances of $11^9$, ramification is tame of order $|I_{11}|=10$. 
 The completions of the two fields at $11$ are actually isomorphic,
 both being $\Q_{11}[\pi]/(\pi^{10}-66) \times \Q_{11} \times \Q_{11}$.  
 This agreement is a little unexpected because $\Q_{11}[\pi]/(\pi^{10}-66)$ is
 one of ten different totally ramified decic $11$-adic fields,
 all with cyclic Galois group $C_{10}$.    
 
 Since the exponents $14$ are at least the degree $12$, ramification
 must be wild at $2$ in $K_1$ and at $3$ in $K_2$.  The database
 describes this ramification in terms of the slope-contents
 $[4/3,4/3]_3^2$ and $[3/2]^2_2$ respectively, the numbers in square brackets
 being wild slopes as explained in \cite{jr-local-database}.  The decomposition
 group at $2$ for $K_1$ is the symmetric group $S_4$, with $A_4$
 being the inertia group, and $V$ being the wild inertia group.  Similarly
 the decomposition group at $3$ for $K_2$ is the dihedral group
 $D_6$, with $S_3$ being the inertia group and $A_3$ the wild 
 inertia group.
 
 It is often enlightening to work not with the discriminant $D$ of a degree $n$ 
 field $K$, but rather with the root discriminant $\delta = |D|^{1/n}$.   
 For example, $\delta$ relates well with the root discriminant $\Delta$ of a 
 Galois closure $L$: one has $\delta \leq \Delta$ with equality
 if and only if $L/K$ is unramified.  For our two cases, the renormalization
 to root discriminants works out to 
 \begin{align}
 \delta_1 & \approx  33.87, &
 \delta_2 & \approx 38.77.
 \end{align}
 The Galois root discriminants are best calculated one prime at a time.  
 If $p$ is tamely ramified with $|I_p|=t$, then its multiplicative contribution is 
 $p^{(t-1)/t}$.  If $t$ is wildly ramified then the contribution can
 be directly computed from the slope content as explained in \cite{jr-local-database}.  In our
 cases, one obtains
 \begin{align*}
 \Delta_1 & = 2^{7/6} 3^{10/11} 11^{9/10} \approx 52.75, &
 \Delta_2 & = 2^{10/11} 3^{7/6} 11^{9/10} \approx 58.55.
 \end{align*}
The root discriminants $\delta_1$, $\delta_2$, $\Delta_1$, and $\Delta_2$
are all unusually small for the Galois group $PGL_2(\F_{11})$, as
we discuss further at various points of Section~\ref{lightly}.
 
 \subsection{Lifts}
 \label{lifts} Serre reciprocity is naturally formulated at the linear level of $GL_2$,
 while we in this paper are working as much as possible at the computationally more
 accessible projective level of $PGL_2$.   To make the connection to modular forms, 
 we first have to lift from the projective level to the linear level.
 
 Let $SL_2^{\pm}(\F_{11})$ be the group of two-by-two matrices over $\F_{11}$ with
 determinant $\pm 1$.  Let $C_5$ be the group of scalar matrices of odd order.  
 Via the product decomposition $GL_2(\F_{11}) = SL_2^{\pm}(\F_{11}) \times C_5$, 
 one can focus attention on $SL_2^{\pm}(\F_{11})$.  This group
 has computational appeal because, unlike $GL_2(\F_{11})$, it is
 a subgroup of $S_{24}$.  
 
For context, note that the polynomials $f_1(x^2)$ and $f_2(x^2)$ both have
Galois group the full wreath product $C_2 \wr PGL_2(\F_{11})$ of order 
$2^{12} \cdot |PGL_2(\F_{11})|$.   At issue is whether the defining polynomials $f_i$
can be adjusted so that replacing $x$ by $x^2$ yields the group
$SL_2^{\pm}(\F_{11})$.  For this to happen, a sign $\epsilon_v \in \{-1,1\}$
has to be $1$ for all the ramified places $\{\infty,2,3,11\}$.   Remarkably
indeed $\epsilon_v=1$ always, and so lifted fields $\tilde{K}_i$ are known to exist.

Finding the better polynomial for the two $\tilde{K}_i$ requires a computation
with $S$-units with $S = \{2,3,11\}$, as discussed in \cite[Chpt~5]{cohen2}.  In the first case,
a polynomial for a lifted field $\tilde{K}_1$ with Galois group $SL_2^{\pm}(\F_{11})$ is  
 \begin{eqnarray*}
\tilde{f}_1(x) \!\!\! & = &  \!\!\! x^{24}-20 x^{22}+208 x^{20}-1380 x^{18}+6432 x^{16}-21696 x^{14} + 52824  x^{12} \\ && -90432 x^{10}+100128 x^8-65728 x^6+31808 x^4-17152 x^2-14256.
\end{eqnarray*}
As $d$ runs over square-free integers, the fields defined by $\tilde{f}_1(\sqrt{d} x)$ run over all lifts of $K_1$.  The field $\tilde{K}_1 = \Q[x]/\tilde{f}_1(x)$
we have chosen has discriminant  $\tilde{D}_1 = -2^{32} 3^{20} 11^{19}$.   Our choice is one of the two with 
smallest discriminant, the other one being given by $f_1(\sqrt{-11} x)$.  
A corresponding $\tilde{f}_2(x) = \tilde{f}_{11}(21,-26,x)$, giving a lifted field discriminant of $\tilde{D}_2 = -2^{20} 3^{30} 11^{19}$, is a specialization
of the parametric family \eqref{tilde11}.

\subsection{Conductors} 
\label{conductors} To make the connection with modular forms, we need to study
the ramification in the lifted fields $\tilde{K}_i$, and then translate to conductors.  
From the discriminants reported above, $\tilde{K}_1/K_1$ is ramified at $2$ and $11$ while
$\tilde{K}_2/K_2$ is ramified at $3$ and $11$.   At the wild prime $2$ in the first case,
an extra wild slope appears so that slope content is now $[3/2,4/3,4/3]_3^2$.  At the
wild prime $3$ in the second case, the tame part of inertia gets larger, so that 
slope content becomes $[3/2]_4^2$.

%
%
%
  The conductors of the 
 Galois representations coming from the inclusion $SL_2^{\pm}(\F_{11}) \subset 
 GL_2(\F_{11})$ are again small:
 \begin{align}
 N_1 & = 24 = 2^3 \cdot 3, &
 N_2 & = 54 = 2 \cdot 3^3.
 \end{align} 
 The source of the exponents $1$ is that $I_p$ in both cases can be taken to be strictly
 upper-triangular matrices, and so the subspace $\F_{11}^2$ fixed by $I_p$ has
 codimension one.    At the primes with exponent $3$, the subspace fixed by $I_p$
 is just $\{0\}$.  However the exponent is $3$ rather than the codimension $2$ because of wildness;
 $3$ arises in both cases as the dimension $2$ times the highest slope $3/2$.  
 

\subsection{Corresponding newforms}  
\label{corresponding}
We will use a standard notation for modular forms, including that
$S_k(N)$ is the space of cusp forms of weight $k$ on the
group $\Gamma_0(N)$.  Via expansion at the cusp $\infty$,
a modular form can be viewed as an element of the power
series ring $\C[[q]]$.  Of particular importance for us
are the new subspaces $S^{\rm new}_k(N)$, which 
has a canonical basis $P_k(N)$ of forms $q + \cdots$ which
are eigenforms for both the Atkin-Lehner operators $w_{p^e}$, 
with $p^e||N$, and the Hecke operators $T_n$, for 
$n \nmid N$.   The word newform always refers to an element
of a $P_k(N)$.  Standard references include \cite{Kob93}, 
\cite{St07}.  Our use of modular forms in this paper is mostly
limited to extracting particular newforms from 
the collection of newforms with rational coefficients drawn 
up in \cite{RobNRC}.  Section~2 of this reference
is a brief synopsis of modular forms, adapted to our current needs.


The Serre reciprocity theorem in our cases says that $K_i$ comes from a newform 
in $P_k(N_i)$, where $k \in \{2,4,6,8,10,12\}$.   Because of the nature of tameness at $11$, Gross's
theory of companion forms of Type 1T says that it comes from two forms, in weights
adding to $12$.  

 In the first case, the sets $P_k(24)$ respectively have size $1$, $1$, $3$, $3$, $5$, and $5$.    
Looking through the eighteen forms, only two match out through $p < 100$, these being
\begin{equation}
\label{first}
{\renewcommand{\arraycolsep}{1.3pt}
\begin{array}{rcrrrrrrrrcc}
g_1 & = & q &+ 3q^3 & + 14q^5 & - 24q^7 & + 9 q^9 & - 28 q^{11} &
 + \cdots & \in & P_4(24), \\
h_1 & = & q & + 27q^3 & - 530q^5 & + 120q^7 & + 729 q^9  & - 7196 q^{11} &
+  \cdots & \in & P_8(24). \\
\end{array}
}
\end{equation}
Here and always for $PGL_2(\F_{11})$ fields, a field $K$ and a newform $\sum a_n q^n \in P_k(N)$ correspond if and only if
 the partition $\lambda_p$ and the normalized square $s_p = a_p^2/p^{k-1} \in \F_{11}$ 
 match for all $p \nmid 11 N$ via  Table~\ref{frobtable}.  
As an example, whenever $\lambda_p = 12$, one must have 
$s_p \in \{7,8\}$.  Similarly, if $s_p=4$, one must have 
$\lambda_p \in \{11^1 1,1^{12}\}$.  As a completely explicit instance
of the correspondence being discussed, consider $K_1$ and $(g_1,h_1)$ at $p=5$.  One has
 $\lambda_5 = 5^2 1^2$ because the irreducible factorization of $f_1(x)$  in $\F_5[x]$
is 
\[
\left(x^5+4 x+2\right) \left(x^5+3 x^3+3 x^2+3 x+2\right) (x+2) (x+4).
\]
Indeed the normalized squares
 $14^2/5^3$ and $530^2/5^7$ both reduce to $5$ in $\F_{11}$, in conformity 
 with Table~\ref{frobtable}.  

The second case is similar.    The sets $P_k(54)$ have sizes $2$, $4$, $6$, $8$, $10$, and $12$.  
Looking through these forms, four match $K_2$ through $p<100$, two of which are
\begin{equation}
\label{gs}
{\renewcommand{\arraycolsep}{3pt}
\begin{array}{rcccrrrrrrrrrrrrcc}
g_2 & =  &q & - q^2&& + q^4 &+ 3q^5 && - q^7 & 
 + \cdots & \in & P_2(54), \\
h_2 & = &  q &+ 16q^2 &&+ 256q^4 & - 435q^5 && - 2527q^7 & 
  + \cdots & \in & P_{10}(54).
\end{array}
}
\end{equation}
The other two which match are twists $g_2^\chi$ and $h_2^\chi$ of the first two,
differing only in that coefficients $a_n$ with $n \equiv 2 (3)$ are negated.  This twisting 
is not seen in the matching criterion.

A subtlety of the general situation is nicely illustrated by looking at weights more closely in 
our pair of examples.   The field $L = \Q_{11}[x]/(x^{10}-66)$ is a splitting field for 
both polynomials $f_1$ and $f_2$.   Let $G_i$ be the Galois group of $f_i$ with respect
to this splitting field, so that both $G_i$ contain the cyclic group $C = \Gal(L/\Q_{11})$ of 
order ten.     Note that there $C$ has four automorphisms $i_j$, where $i_j(\sigma) = \sigma^j$, 
and $j$ can be $1$, $3$, $7$ or $9$.  The normalizer of $C$ is $G_i$ is $D_i$, a dihedral group of
order $20$.   This means that of the $|PGL_2(\F_{11})| = 12 \cdot 11 \cdot 10$ isomorphisms
from $G_1$ to $G_2$, twenty take $C$ to $C$.  Ten of these are some $i_j$ and 
ten are $i_{j'}$, with $j+j'=10$.   A Galois-theoretic computation
shows that $\{j,j'\} = \{3,7\}$, not the other other possibility, $\{1,9\}$.  
This is why the two weight sets $\{4,8\}$ and $\{2,10\}$ are different.    
Refining this computation further, one can see purely Galois-theoretically that 
the weights for $f_1$ and $f_2$ have $\{4,8\}$ and $\{2,10\}$ respectively.

\section{Fifty-three $PGL_2(\F_\ell)$ fields with rational companion forms}
\label{fifty-three}

    In this section, we prove the existence of fifty-three $PGL_2(\F_\ell)$ fields $K$
with associated rational companion forms $g$ and $h$.   In contrast to the previous section,
here we start with $(g,h)$ and obtain only the abstract existence of $K$, not an 
explicit defining polynomial $f(x)$.  The examples of the previous section
provide helpful illustrations, but to a large extent this section
can be read independently.


\subsection{Triples $(K,g,h)$}
\label{triples}
    In the introduction, we explained that we are seeking fields $K$ belonging to triples $(K,g,h)$
similar to the triples $(K_i,g_i,h_i)$ of the previous section.   In this subsection, we define precisely the
type of triples we seek.  

    Throughout, a prime $\ell$ is present, typically not incorporated into the notation. 
     We exclude the prime $\ell=2$ because it 
behaves slightly differently.  When we pursue classification starting in \S\ref{obtaining},
 we will take $\ell \geq 11$.    The fields $K$ we allow
are those presentable of the form $\Q[x]/f(x)$ with $f(x) \in \Q[x]$ a degree $\ell+1$
 polynomial with Galois
group $PGL_2(\F_\ell)$.  

    Beyond the prime $\ell$, three more invariants associated to a triple $(K,g,h)$ are positive
integers $N$, $k$, and $k'$.   The integer $N$, 
called the level or the conductor, is required to be not a multiple
of $\ell$.  The integers  $k$ and $k'$, called the weights, are
 even and in the range $[2,\ell+1]$.
The remaining entries of the triples
are newforms with the common level $N$, trivial character, rational 
coefficients, and the indicated weights:
\[
g = \sum a_n q^n \in S^{\rm new}_k(N), \;\;\;\;\;\;\;\;\;
h = \sum b_n q^n \in S^{\rm new}_{k'}(N).
\]
Here basic notation for newforms has been recalled in a formalistic 
way at the beginning of \S\ref{corresponding}, with references
also given there.   

     A triple $(K,g,h)$ has yet more invariants.
In particular, for primes $p$ not dividing $N$, the field
$K$ yields a factorization partition $\lambda_p$, 
the form $g$ determines a normalized square $s_p = a^2_p/p^{k-1} \in \F_\ell$, 
and the form $h$ likewise determines a normalized square $s'_p = b^2_p/p^{k'-1} \in \F_\ell$.  These numeric
invariants are required to correspond as follows.    First, $s_p = s_p'$.  Second,
let $o_p$ be the least common multiple of the parts of $\lambda_p$.  Let
$O_p$ be the common order of the elements 
\[
\left(
\begin{array}{cc}
0 & - 1 \\
p^{k-1} & a_p
\end{array}
\right) 
\mbox{ and } 
\left(
\begin{array}{cc}
0 & - 1 \\
p^{k'-1} & b_p
\end{array}
\right) 
\]
in the group $PGL_2(\F_\ell)$.  Then it is required that either $o_p = O_p$ 
or $(o_p,O_p)=(1,\ell)$.   Table~\ref{frobtable} illustrates this correspondence
for $\ell=11$. 

We require that $g$ and $h$ are related via cyclotomic twisting as follows.  
For either $t=1$ or $t=2$, we require that $k+k' = \ell-1+2t$ and  
\begin{equation}
\label{congruence}
n^{t} a_n \equiv n^k b_n \; (\ell)
\end{equation}
for all $n$.  Symmetry between $g$ and $h$ is present, because this last condition could be equivalently rewritten as  
$n^{k'} a_n \equiv n^{t} b_n \; (\ell)$.   If $k = k'$, either congruence says that 
$a_n \equiv \chi(n) n$, with $\chi(\cdot) = (\cdot/\ell)$ 
the quadratic character on $\F_\ell$.    As our last requirement on triples, we partially normalize by demanding
$k \leq k'$. 

    Given a triple $(K,g,h)$, one can sometimes trivially make a new one in three
ways.  First, one can replace $g$ and/or $h$ by a different form with the same
reduction to $\F_\ell[[q]]$ in the somewhat rare case that such a form exists.  
Second, one may be able to twist $g$ and/or $h$ by a quadratic character,
keeping the level the same, as discussed with examples after \eqref{gs}. 
Finally if $k=k'$ one can simply switch $g$ and $h$.  Note that in this
final case, $g$ and $h$ are not trivially obtained from one another because
twisting either $g$ or $h$ by the quadratic character $\chi$ above 
increases the level from $N$ to $\ell^2 N$.  
   None of these operations change $K$, which
is why Theorem~\ref{mt} below counts $K$ extendable to triples,
rather than triples themselves.

  \subsection{Conditions on $K$ and companion forms} 
  \label{conditions}
This subsection discusses which $PGL_2(\F_\ell)$ number fields $K$ have a chance of being 
in a triple $(K,g,h)$.   Our discussion defines the three types 
of companion forms, and in the process motivates some of the
definitions made in the previous subsection.
The next four paragraphs 
do not use rationality.  

Let $K$ be any $PGL_2(\F_\ell)$ field which is not totally real.  
Then, by Serre reciprocity,  $K$ comes from a newform $g$, perhaps with irrational coefficients,
in some space $S^{\rm new}_k(N,\chi)$.  Here the Dirichlet character can be non-trivial in 
general, and the weight satisfies $\chi(-1) = (-1)^k$ and is therefore allowed to be odd.  Matching
is between $\lambda_p$ as before and now $s_p = a_p^2/(p^{k-1} \chi(p))$.     One
can always take $k$ in the interval $[2,\ell+1]$, which accounts
for our making this restriction on $k$ in the previous section.    

Let $D$ be the discriminant of $K$, and write $c = \ord_{\ell}(D)$.   The largest
$c$ can be is $2 \ell-1$. Suppose $c>\ell$, so that $c$ has the form $\ell-2+k$ for $k \in [3,\ell+1]$.
Then $g$ necessarily has weight $k$.  The condition \eqref{congruence} makes sense for 
irrational forms reduced to characteristic $\ell$ as well, and there is no companion form $h$.  
Examples with an explicit polynomial $f(x)$ for $K$ and $g$ rational are given 
for $N=1$ and $\ell \geq 11$  in \cite{Bos11}, and for 
$N \in \{2,3,4,6,8\}$ and $\ell \leq 7$ in \cite[\S6.3]{RobNRC}.  

Still allowing irrational forms and general character, 
 consider the complementary case $c \leq \ell$,  except 
that we temporarily exclude $c=0$.  For these more lightly ramified 
fields, there is always a companion form $h$ with weight $k'$ satisfying $k+k'=\ell-1+2t$ with
$t \in \{1,2\}$ as for \eqref{congruence}.  
The three cases are as follows:
  \begin{equation}
  \label{conditionstab}
  \begin{array}{l |ccc}
  \mbox{Case} & \mbox{Size of $\ell$-inertia} & \ord_\ell(D) &  \mbox{Consequence for $(k,k')$} \\
  \hline
  \mbox{1T=diagonalizable} & (\ell-1)/d & \ell-1-d & k+k' = \ell+1\\
   \mbox{2T=supersingular} & (\ell+1)/d &  \ell+1-d & k+k' = \ell+3 \\
  \mbox{2W=peu ramifi\'ee} & \ell(\ell-1) & \ell  & (k,k') = (2,\ell+1).  
  \end{array}
  \end{equation}
  The case $c=0$ would be similar, except that the two weights should be
  $(k,k') = (1,\ell)$, with the form of weight $1$ usually living only in 
  characteristic $\ell$ \cite{CV92}.  
  Summarizing, a fundamental reason to be interested in
  companion forms from a number-theoretic viewpoint is 
  that existence of a companion form translates into light ramification at $\ell$.  
  
  An important distinction between the cases 
  deserves to be mentioned.   If from $g$ one sees that $s_\ell$ is zero in $\F_\ell$,
  then one is automatically in Case~2T; one knows without looking that there
  is then a companion form $h$ in weight $\ell+3-k$.     However if $s_\ell \neq 0$,
  then one really needs to look for a companion form 
  in weight $\ell+1-k$ to identify ramification.  One is in Case~1T only if there is a such a companion
  form.  Otherwise, $c = \ell-2+k$ as above, and ramification is 
  wild.  If $k=2$, one is in the case 2W and there is necessarily 
  a companion form, but now in weight $\ell+1$ rather than $\ell-1$.  If $k>2$, 
  one is in the generic case 
  and there is no companion form.   
  
  Now we return to our requirement that both $g$ and $h$ have rational
  coefficients.  Then the Nebentypus $\chi$ is trivial, forcing the discriminant $D$ to be
  $(-1)^{(\ell-1)/2} \ell$ times a square.    Moreover the field 
  $K$ has a lift to a field $\tilde{K}$ with Galois group embedding in
  $GL_2(\F_\ell)$.    For any $\ell$ we expect infinitely many $K$ 
  to satisfy these two natural conditions.  They all fit into a $(K,g,h)$, as long
  as we allow irrational $(g,h)$.    Besides implying these two conditions,
  rationality of coefficients on the modular side does 
  not translate into anything natural on the Galois side. 
  Rather it just corresponds to restricting to a presumably
  much smaller subset of this collection of fields $K$.

\subsection{Finding and confirming triples $(K,g,h)$}
\label{obtaining}
To obtain triples $(K,g,h)$ with $g$ and $h$ rational, we use the collection of rational newforms without complex multiplication 
built up in \cite{RobNRC}.   
\begin{table}[htb]
{
\[
{\renewcommand{\arraycolsep}{3.8pt}
\begin{array}{r|ccccccccccccc}
\ell & k=2 & 4 & 6 & 8 & 10 & 12 & \;14\; & \;16\; & \;18\; & \;20\; & \;22\; \\
 \hline
11 &  \begin{array}{c} { {\bf 54}} \\ 182 \end{array} &  \multicolumn{1}{c|}{\begin{array}{c} {{\bf 24}}\\ 42 \\ 120
 \end{array}}&  \multicolumn{1}{c|}{\cellcolor{lg} \begin{array}{c} 78 \\ 78 \end{array}}  &  \begin{array}{c} {{\bf 24}} \\ 42\\ 
120 \end{array} &  \multicolumn{1}{c|}{ \begin{array}{c} {{\bf 54}} \\ 182 \end{array} } & \text{} & \text{} & \text{} & \text{} & \text{} & \text{} \\
\cline{4-4} \cline{7-7}
13 &   & 22 & \multicolumn{1}{c|}{\begin{array}{c} 7 \\ 70 \\ 84 \\ 210 \end{array}} & \begin{array}{c} 7 \\ 70 \\ 84 \\ 210 \end{array} & {22}  & \multicolumn{1}{c|}{}& \text{} & \text{} & \text{} & \text{} & \text{}  \\
\cline{5-5} \cline{8-9}
17 &   &  & 24 &  \multicolumn{1}{c|}{42} & 42 & 24 & &  \multicolumn{1}{c|}{} & \text{} & \text{} &   \\
\cline{6-6} \cline{10-10}
19 &   & 10 &  &  \multicolumn{1}{c|}{24} &  \multicolumn{1}{c|}{\cellcolor{lg}}& 24 & \text{} & 10 &  \multicolumn{1}{c|}{} & \text{} & \text{}  \\
\cline{6-7} \cline{11-12}
23 &   &  &  &  &  \multicolumn{1}{c|}{30}&  \multicolumn{1}{c|}{\cellcolor{lg}}& 30 & \text{} & \text{} & \text{} & \multicolumn{1}{c|}{}  \\
\cline{7-8}
29 &  &  &  &  &  & 6 &  \multicolumn{1}{c|}{8} & 8 & 6 & \text{} & \text{}  \\
\cline{9-9}
31 &   &  &  &  &  &  &  \multicolumn{1}{c|}{12} &  \multicolumn{1}{c|}{\cellcolor{lg}} & 12 & \text{} & \text{}  \\
\cline{9-10}
37 &   &  &  &  \;\;\;\;\;\;\;\;\; &  \;\;\;\;\;\;\;\;\; &  \;\;\;\;\;\;\;\;\; &  \;\;\;\;\;\;\;\;\; & \;\;\;\;\;\;\;\;\;  &  \multicolumn{1}{c|}{\;\;\;\;\;\;\;\;\; }  &  \;\;\;\;\;\;\;\;\; &  \;\;\;\;\;\;\;\;\;  \\
\cline{11-11}
41 & \;\;\;\;\;\;\;\;\; &  \;\;\;\;\;\;\;\;\;&  \;\;\;\;\;\;\;\;\;   &  &  &  &  &  &  & \multicolumn{1}{c|}{3} & 3 \\
\cline{12-12}
\end{array}
}
\]
}
\caption{\label{Type1} Guide to the twenty triples $(K,g,h)$ of Type 1T.  The boldface
{\bf 24}'s and {\bf 54}'s represent
the triples of Section~\ref{two}.}
\end{table}
 For $t \in {1,2}$ and each $(\ell,N,k,k')$ within the range of the collection, we 
look at all potentially congruent rational pairs $g \in S_k^{\rm new}(N)$ and 
$h \in S_{k'}^{\rm new}(N)$. 
\begin{table}[htb]
{
\[
{\renewcommand{\arraycolsep}{3.9pt}
\begin{array}{c|cccccccccccccccccccccccc}
\ell & k=2& 4 & 6 & 8 & 10 & 12 & \;14\; & \;16\; & \;18\; & \;20\; & \;22\; & \;24 \;   \\
\cline{1-13}
  11 &  \begin{array}{c} 14 \\ w15 \\ 20 \\ w24 \\ 30\\ w42 \\ w84 \\ w96 \end{array} 
   &  \begin{array}{c} (8) \end{array}   & \multicolumn{1}{c|}{  \begin{array}{c} 10 \\ 42 \\ 70 \\ 
  78\\ 96 \end{array} } &  
   \begin{array}{c} 10 \\ 42 \\ 70 \\ 78 \\ 96 \end{array}  &  \begin{array}{c} (8) \end{array}    &
 \multicolumn{1}{c|}{  \begin{array}{c} 14 \\ w15 \\ 20 \\ w24 \\ 30 \\ w42 \\ w84 \\ w96 \end{array}  }& \text{} & \text{} & \text{} & \text{} & \text{} &
   \text{} & \text{} & \text{} & \text{} & \text{} & \text{} & \text{} & \text{} &
  \text{} & \text{} & \text{}  \\
  \cline{8-8}
13 &  \text{} & \text{} &  \begin{array}{c} \mathbf{5} \\ 24 \\ 38 \\ 50 \\ 54 \\ 294 \end{array}  &\cellcolor{lg} \begin{array}{c} 30 \\ 30 \end{array} &    \begin{array}{c} \mathbf{5} \\ 24\\ 38
\\ 50 \\ 54 \\ 294 \end{array}  & \text{} &
 \multicolumn{1}{c|}{ \text{}} & \text{} & \text{} & \text{} & \text{} & \text{} & \text{} & \text{} &
   \text{} & \text{} & \text{} & \text{} & \text{} & \text{} & \text{} & \text{} 
  \\
   \cline{9-10}
17 &   w30  & \text{} &  30  &  \begin{array}{c} 50 \\ 66 \end{array} & \cellcolor{lg} \begin{array}{c} 42\\ 42 \end{array}
 &   \begin{array}{c} 50 \\ 66 \end{array} &  30  & \text{} &   \multicolumn{1}{c|}{w30}
   & \text{} & \text{} & \text{} & \text{} & \text{} & \text{} & \text{} & \text{} &
   \text{} & \text{} & \text{} & \text{} & \text{} \\
   \cline{11-11}
19 &  \text{} &  30  &  \begin{array}{c} (4) \end{array}   & \text{} &  \multicolumn{1}{c|}{ \begin{array}{c} \mathbf{3} \\ 90 \end{array}}   &  \begin{array}{c} \mathbf{3} \\ 90 \end{array}   & \text{} &  \begin{array}{c} (4) \end{array}   &  30  &
   \multicolumn{1}{c|}{\text{}} & \text{} & \text{} & \text{} & \text{} & \text{} & \text{} & \text{} &
   \text{} & \text{} & \text{} & \text{} & \text{}  \\
   \cline{7-7}  \cline{12-13}
23 &  \text{} & \text{} & \text{} & \text{} & \text{} &  \multicolumn{1}{c|}{\text{}} &  \text{} & \text{} & \text{}
   & \text{} & \text{} &  \multicolumn{1}{c|}{\text{}}& \text{} & \text{} & \text{} & \text{} & \text{} &
   \text{} & \text{} & \text{} & \text{} & \text{}  \\
   \cline{8-8} 
 29 & \text{} & \text{} & \text{} &  30  &  12  & \text{} &  \mathbf{2}  & \cellcolor{lg}  \text{} &  \mathbf{2} &
   \text{} &  12  &  30  & \text{} & \text{} & \text{} & \text{} & \text{} & \text{}
   & \text{} & \text{} & \text{} & \text{}  \\
31 &  \text{} & \text{} & \text{} & \text{} & \text{} & \text{} & \text{} &  \multicolumn{1}{c|}{\text{}} & \text{} & \text{}
   & \text{} & \text{} & \text{} & \text{} & \text{} & \text{} & \text{} &
   \text{} & \text{} & \text{} & \text{} & \text{}  \\
   \cline{10-10}
37 &  \text{} & \text{} & \text{} & \text{} & \text{} &  6  & \text{} & \text{} & \text{} &
    \cellcolor{lg}  \text{} & \text{} & \text{} & \text{} & \text{} & \text{} & \text{}
   & \text{} & \text{} & \text{} & \text{}  \\
41 &  \text{} & \text{} & \text{} & \text{} & \text{} &    & \text{} & \text{} & \text{} &  \text{} &
    \cellcolor{lg}  \text{} & \text{} & \text{} &    & \text{} & \text{} & \text{} & \text{}
   & \text{} & \text{} & \text{} & \text{}  \\
 43 &  \text{} & \text{} & \text{} & \text{} & 6 &   & \text{} & \text{} & \text{} & \text{} &
     \multicolumn{1}{c|}{\text{}} &  \text{} & \text{} &   & \text{} & \text{} & \text{} & \text{}
   & \text{} & \text{} & \text{} & \text{}  \\

\end{array}
}
\]
}
\caption{\label{Type2} Guide to the twenty-seven triples of type 2T and the six triples 
of type 2W.   To the right of the chart, there should also be $6$'s at locations $(\ell,k) = (37,28)$ and $(43,36)$
 The boldface {\bf 5}'s, {\bf 3}'s, and {\bf 2}'s are pursued in \S\ref{four}.  
Also in parentheses are two degenerate triples discussed in \S\ref{twoquartics}}
\end{table}
 Motivated by the definitional congruences \eqref{congruence}, we consider
$\delta = \sum_n c_n q^n \in \F_{\ell}[[q]]$
with $c_n = n^t a_n - b_n n^{k}$.  
Let $\theta = q \frac{d}{dq}$ be Ramanujan's theta operator.  As described in \cite[(4.5)]{gr}, this operator increases
weights of reduced modular forms in  $\F_\ell[[q]]$ by $\ell+1$.  The difference $\delta = \theta^2 f - \theta^{k+2} g$ in question is 
then the reduction of a modular form in $M_\kappa(N)$, where $\kappa = k' + (k+t)(\ell+1)$.   
Let $\sigma_1(N) = \prod_{p^e||N} (p^e+p^{e-1})$ be the index of $\Gamma_0(N)$.  Then 
by \cite[Thm~9.18]{St07}, $\delta$ is determined by its Fourier coefficients
$c_n$ for $n$ at most the Sturm bound $S = \kappa \sigma_1(N)/12$.    We compute these $c_n$ until either one of them is nonzero or 
$S$ is reached.  In the latter case, we have confirmed that indeed $(g,h)$ is a companion pair.   Sometimes
the common  $s_p = a_p^2/p^{k-1} = b_p^2/p^{k-1}$ inspected for $p \nmid \ell N$ do 
not suffice to ensure the surjectivity of the projective representation.  In these few cases, we 
identify the relevant number field, so as to unconditionally confirm lack of 
surjectivity.   Numerics associated to the fifty-three pairs remaining are in 
Tables~\ref{Type1} and \ref{Type2}.  Each field $K$ gives rise 
to a minimal conductor $N$ appearing twice
in the $\ell$-row, once in the $k$ column and once in the $k'$ column.      
The parenthesized entries $(N)$
on Table~\ref{Type2} indicate two of the $(g,h)$ discarded because of
nonsurjectivity of the Galois representation into $PGL_2(\F_\ell)$.   

As an example, the largest Sturm bound encountered in a Type 1 pair 
occurs for $(\ell,k,k',N) = (13,6,8,210)$.  Here the forms are the unique newforms in their
one-dimensional Atkin-Lehner eigenspaces, $g \in S^{\rm new}_{6}(210)^{--++}$
and $h \in S^{\rm new}_{8}(210)^{+---}$.  One has $\kappa = 92$, $\sigma_1(210)=576$, and so the Sturm bound is $S=4416$.
A  {\em Magma} computation using built-in commands in a straightforward 
way took $24$ minutes to confirm that $\delta \in \F_{13}[[q]]$ is indeed zero.

Summarizing, our computations prove the following theorem. 
\begin{Theorem}  \label{mt}  There are at least fifty-three number fields $K$ belonging to triples
$(K,g,h)$ satisfying the following conditions: $K$ has degree $\ell+1$ and
associated Galois group $PGL_2(\F_\ell)$ for a prime $\ell \geq 11$;  
$g$ and $h$ are rational newforms which are companion forms
modulo $\ell$ and have projective modulo $\ell$ representations
corresponding to $K$.  
\end{Theorem}
In the sequel, we will sometimes simply speak of the fifty-three triples $(K,g,h)$,
always chosen with $g$ and $h$ having minimal level $N$, even 
though in a few cases there are the ambiguities in $g$ or $h$ mentioned
at the end of \S\ref{triples}.  
%
\subsection{Conjectural completeness}    We believe that the
list of fifty-three number fields in Theorem~\ref{mt} is complete.  In this subsection,
we give our reasons; in brief, the fifty-three number fields arise
towards the the beginning of our search.  We have computed much
further and found no more.    

In general, for the weights $k \leq k'$ associated to the pair $(g,h)$, one has 
$k' \geq (\ell+1)/2$.   Conjecture~1.1 of \cite{RobNRC} says that there are 
no non-CM newforms with
rational coefficients 
and weight $k' \geq 52$.    This would imply that there are indeed no fields 
for $\ell \geq 101$.   Conjecture~1.1 of \cite{RobNRC}
 says moreover that all such newforms of weight $k' \geq 18$ are known,
as indeed their minimal levels $N$ are always $\leq 30$.  This would imply our
list of three fields for $\ell \geq 37$ is complete.    It is moreover argued in \cite{RobNRC} that 
all or very close to all such newforms in weights $10$, $12$, $14$, $16$ are known too.  
This makes our list twenty-one fields for $\ell \geq 17$ likely to be complete too. 

The evidence in \cite{RobNRC} suggests that for $k \in \{6,8\}$, 
there are likely a few non-CM rational newforms with minimal level 
beyond the cutoffs $C_6 = 1000$ and $C_8 = 700$ used there.  
However it seems unlikely to us that these unknown newforms 
are part of a companion pair, expecially given that 
the largest level $N$ appearing on Tables~\ref{Type1} and \ref{Type2} 
is $294$.  It is for this reason that we believe that our
lists of thirteen and nineteen fields for $\ell = 13$ and $\ell = 11$ respectively 
are complete
as well.

\subsection{Explicit formulas} In \cite{RobNRC}, we explained how to 
get completely explicit formulas for newforms in the cases 
$N \in \{2,3,4,6,8\}$.   The entries for these $N$ in Tables~\ref{Type1} and \ref{Type2} 
correspond to the following companion forms
\begin{align*}
N=2:  &&  \Delta^+_{14} & \stackrel{29}{\sim}  \Delta^-_{18}  \\
N=3:  && \Delta^-_{10} & \stackrel{19}{\sim}  \Delta^+_{12}, &  \Delta^-_{20} & \stackrel{41}{\sim}  \Delta^{+a}_{22},   \\
N=4:  &&  (\Delta^-_{6} & \stackrel{19}{\sim}  \Delta^-_{16}),\\
N=6:  &&  \Delta^{--}_{12} & \stackrel{29}{\sim}  \Delta^{++}_{18},  &
   \Delta^{--}_{12} & \stackrel{37}{\sim}  \Delta^{+-}_{28},  &  \Delta^{+-}_{10} & \stackrel{43}{\sim}  \Delta^{-+}_{36}, \\
N=8:  &&  (\Delta_{4}^+ & \stackrel{11}{\sim}  \Delta_{10}^-), &  \Delta_{14}^+ & \stackrel{29}{\sim}  \Delta_{16}^{-b}.
\end{align*}
Here $\Delta_k^\epsilon$ denotes the unique newform of weight $k$ on $\Gamma_0(N)$ with Atkin-Lehner eigenvalue 
string $\epsilon$.  The unusual cases $(k,N) = (3,22)$ and $(8,16)$ are highlighted in
\cite{RobNRC}, as in these cases a two-dimensional Atkin-Lehner space has two rational
newforms.  The above display identifies which newform is involved in the compansion 
forms, with completely explicit formulas for $\Delta^{+a}_{22}$ and $\Delta^{-b}_{16}$ 
given in \S5.3 and \S5.6 of \cite{RobNRC} respectively.

One could also write down explicit formulas for other $N$.    For example, define $\Theta_t = \sum_{x,y \in \Z^2} q^{t(x^2+xy+y^2)}$.
Returning to our very first example with $N=24$, the unique newform of weight two is 
$\Delta_2^{-+} = 18^{-1} (\Theta_4-\Theta_1)(\Theta_2 - 4\Theta_8)$.  Like all the other generators of cuspidal ideals
considered in \cite{RobNRC}, it is also an eta-product, $\eta_2 \eta_4 \eta_6 \eta_{12}$.  The companion
forms from \eqref{first} have the explicit formulas
\begin{align*}
g_1 & = \Delta_{4}^{--}  = 3^{-1} (\Theta_2^2 + 2\Theta_4^2)  \Delta_{2}^{-+},  \\
 h_1 & = \Delta_{8}^{+-} = 9^{-1} (\Theta_2^2 - 2\Theta_4^2)(7\Theta_2^4 - 44\Theta_2^2\Theta_4^2 + 28\Theta_4^4) \Delta_{2}^{-+}.
\end{align*}
The formulas constructible from \cite{RobNRC} for the nine companion pairs displayed above 
are of a similar nature.  The main congruence \eqref{congruence} can be seen explicitly, by expanding the power series.  
For example, for both
$g_1$ and $h_1$, the first five primes having Fourier coefficient congruent to zero modulo $11$
are 103, 149, 179, 197, and 257.

\section{Lightly ramified number fields}   
\label{lightly}
    This section first obtains polynomials for some of the fifty-three number fields 
 from the last section.  It next analyzes 
ramification in these number fields, finding that some root discriminants are particularly small.    
Finally it discusses the natural problem of obtaining complete lists of number
fields with Galois group a finite subquotient of $GL_2(\overline{\F}_\ell)$ and small root discriminant.



\subsection{Polynomials from elliptic curves} 
\label{Polynomials1}
For eleven of our triples $(K,g,h)$,
the modular weight of the form $g$ is $2$.   As $g$ has rational coefficients, there is a corresponding
elliptic curve $E_g$, easily found on the LMFDB \cite{LMFDB}.  
 A degree $\ell+1$  polynomial can then be obtained for $K$ 
by looking at the $\ell+1$ different subgroups of $E_g(\C)$ of size $\ell$.  
{\em Magma}'s \verb@AtkinModularPolynomial@ does this immediately.  

For example, only one of eleven triples has $\ell \neq 11$,
and its residual prime is $\ell = 17$.     There are eight elliptic curves with the required
conductor $30$, all isogenous and hence all yielding the correct $K$.  One of the
$j$-invariants is $71^3/(2^4 3^3 5)$.   Specializing the Atkin modular polynomial to
this $j$-invariant and applying \verb@polredabs@ yields
\begin{eqnarray*}
&& x^{18}-2 x^{17}-714 x^{15}+3060 x^{14}-7854 x^{13}+258468 x^{12}-1062840
    x^{11} \\ && -2425764 x^{10}-6360720 x^9+224396532 x^8-694308084 x^7 \\ && -1149382920
    x^6+8831732832 x^5-3417673200 x^4-68962605552 x^3 \\ && +206896699224
    x^2-267387716040 x+143842600848,
\end{eqnarray*}
with discriminant $2^{16} 3^{16} 5^{16} 17^{17}$.  
The large coefficients are a reflection of the relatively large root discriminant  
$\delta \approx 298.6$. 

The other ten triples all have $\ell=11$ and so can be treated uniformly,
even at the lifted level of degree $24$ fields $\tilde{K}$ with Galois 
group $SL^{\pm}_2(\F_{11})$.  Here it does not suffice to work with
 $j$-invariants, as quadratic twisting is seen in the lift.  Accordingly we work with
 actual elliptic curves $y^2 = x^3 + a x + b$.   Starting
 with the Atkin modular polynomial with $\ell=11$, lifting
 to degree $24$ polynomials for individual $j$, and interpolating, 
 we obtain the following polynomial with just seven terms:
   \begin{eqnarray}
\nonumber  \tilde{f}_{11}(a,b,x) & = & d^{10}  x^{24} -15840 d^5 x^{12} -337920 a d^3 x^8 \\ && 
\label{tilde11} \qquad    -2280960 b d^2 x^6 +    811008 a^2 d x^4+663552 a b x^2-2816 .
    \end{eqnarray}
 Here we have abbreviated using the discriminant $d=-4 a^3 - 27 b^2$.  
 Correctness of the seven-term polynomial is algebraically confirmed by  
 comparing with a full 
 $11$-division polynomial of degree $120$ and factoring a resolvent.   The polynomial
 applies to our ten cases through the following chart:
 \[
 \begin{array}{c|rr|rrr|rrrrr}
   & \multicolumn{2}{c|}{\mbox{1T}} &  \multicolumn{3}{c|}{\mbox{2T}} &  \multicolumn{5}{c}{\mbox{2W}} \\
N  & 54 & 182 & 14 & 20 & 30 &  15 & 24 & 42 & 84 & 96   \\
\hline
a   & 21 & 13861   & -675  &  -108 & 1917 & -27 & 54 & - 5211 &-108 &  -189  \\
b   & -26 & 426358  & 13662 & 297 &  99198 & 8694 & 189 & 319734 & -1755 & -540 
\end{array}.
\]
We remark that we have found that Atkin modular polynomials for several other $\ell$ also have 
fewnomial equivalents;  in the cases $\ell \equiv 3 \; (4)$,  some of these lift to 
$SL^{\pm}_2(\F_\ell)$-covers via $x \mapsto x^2$. 

\subsection{Polynomials from higher weight modular forms}  
\label{polynomials2}
For forty-two of our $(K,g,h)$, the smaller weight $k$ is at least four.    Our discussion
so far has included a polynomial for only one of these fields $K$, namely our very first example
$K_1$, with defining polynomial  $f_1(x)$ from \eqref{examp1}.   

It is however theoretically possible to take a modular form as a starting point and
compute an associated mod $\ell$ Galois representation.  
Explicit examples in the literature currently start from either forms with rational coefficients in 
level $N=1$ \cite{bosman-level1,Mas13} 
or forms with irrational coefficients in weight $k=2$ 
\cite{bosman-sl216, bosman-inverse}.  

Our collection of examples provides a testing ground for these methods in the setting of
$N>1$ and $k>2$.   Most of them seem to be currently beyond computational reach.
However Mascot has very recently computed two new polynomials
for $PGL_2(\F_{13})$, one from the $N=7$ entry on Table~\ref{Type1} and 
one from the $N=5$ entry on Tables~\ref{Type2} and~\ref{fourtable}.   The computation
is explained in detail in \cite{Mas16}, and passes through explicit degree $56$ polynomials
giving quartic lifts.




\subsection{Ramification in modular fields} 
\label{ramificationmod}
One does not actually need polynomials to determine 
ramification in our fields $K$, as Serre reciprocity
is refined enough to calculate it on the modular side.  

For a prime $p \neq \ell$ exactly dividing a minimal conductor $N$, ramification
is tame of order $\ell$.  It contributes 
$p^{(\ell-1)/(\ell+1)}$ to the root discriminant and $p^{(\ell-1)/\ell}$ towards
the Galois root discriminant. If $p^2$ exactly divides the minimal
conductor $N$, then ramification is tame of order $e = 3$, $4$, or $6$,
so that $e$ divides exactly one of $\ell-1$ or $\ell+1$.   If $p=2$ the only
possibility is $e=3$ and if $p=3$ the only possibility is $e=4$, as 
otherwise ramification would be wild.   The contribution to the root discriminant
 is $p^{(e-1)(\ell-1)/(e(\ell+1))}$ if 
$e$ divides $\ell-1$ and $p^{(e-1)/e}$ if $e$ divides $\ell+1$.  
The computation to the Galois root discriminant is always $p^{(e-1)/e}$.   
If $\ord_p(N) \geq 3$, ramification is wild and the procedure is 
more complicated, as with the two examples given in Section~\ref{two}.  

The contribution of $\ell$ to the root discriminant depends on the 
type and weights.   Taking \eqref{conditionstab} as a starting point, 
and writing $-$ in type 1T and $+$ in type 2T, put $d=\gcd(k-1,\ell \pm 1)=\gcd(k'-1,\ell \pm1)$ and $e = (\ell \pm 1)/d$.  
The size of the inertia group $I_\ell$ is then $e$.   The contributions of $\ell$ 
are $\ell^{(e-1)(\ell \pm 1)/(e(\ell+1))}$ 
to the root discriminant and $\ell^{(e-1)/e}$ to the Galois root discriminant. 
As $k$ and $k'$ are even, $d$ is always odd.     
In fact $d=1$ except for the cases $(\ell,N)=(11,78), (13,22), (19,10)$ 
on Table~\ref{Type1},  $(13,30)$, $(17,42)$, $(29,12)$ on Table~\ref{Type2}, 
and the degenerate cases $(11,8)$ and $(19,4)$ on Table~\ref{Type2}.   
Here the inertial group size reductions are respectively $d = 5$, $3$, $3$, $7$, $9$, $3$, $3$, and $5$.   For cases of type 2W, the contribution to the root discriminant is
$\ell^{\ell/(\ell+1)}$, while the contribution to the Galois root discriminant is
$\ell^{1+(\ell-2)/(\ell(\ell-1))}$.  

\subsection{Four lightly ramified fields}
\label{four}
The root discriminant $\delta$ and Galois root discriminant $\Delta$ for four of our fields are
in the middle block of Table~\ref{fourtable}. 
\begin{table}[htb]
\[
\begin{array}{c|c  rr | rrr }
      &   \multicolumn{3}{c|}{\mbox{Tame at $\ell$ from}} & \multicolumn{3}{c}{\mbox{wild at $\ell$ from } }\\
            &   \multicolumn{3}{c|}{\mbox{ companion forms}} & \multicolumn{3}{c}{\mbox{Ramanujan's $\Delta_{12}$} }\\
\ell & N & \delta  & \Delta & N & \delta & \Delta  \\
\hline
11 & 24 & 33.87  & 52.75& 1& 66.44 &  118.39  \\
13 & 5 &  {43.00} & \mathit{47.82}  &1&67.62& 112.04 \\
19 & 3 &  \mathit{44.07} & \mathit{46.43}  &1&71.48& 103.60  \\
29 & 2 &  \mathit{49.50} & \mathit{50.62} &1&79.64&  103.59  
\end{array}
\]
\caption{\label{fourtable} Root discriminants $\delta$ and Galois root discriminants $\Delta$ for 
eight lightly ramified number fields with Galois group $PGL_2(\F_\ell)$.  Italicized entries
are candidates are for smallest possible for their context.  } 
\end{table}
The last three cases are uniformly behaved as they all have type 2T with $N=p$ prime.
Their root discriminants are given by $p^{(\ell-1)/(\ell+1)} \ell^{\ell/(\ell+1)}$ while
their Galois root discriminants are given by the slightly larger number 
$p^{(\ell-1)/{\ell}} \ell^{\ell/(\ell+1)}$.  The rest of this subsection puts the 
four pairs 
$(\delta,\Delta)$ into context.  

\subsubsection{Comparison with the Serre-Odlyzko constant $\Omega$}
Analytic lower bounds on root discriminants of degree $n$ fields increase as
$n \rightarrow \infty$ to an asymptotic limit of $\Omega' = 4 \pi e^\gamma$, with
$\gamma \approx 0.5772$ being Euler's gamma constant.  Under the generalized Riemann
hypothesis, these bounds are increased so that the asymptotic limit becomes 
$\Omega = 2 \Omega' \approx 44.76$ \cite[(2.6)]{od-survey}.  In \cite{jrlowgrd} and then \cite[\S9,10]{jr-global-database}, we put 
forward the principle that it is extremely unusual for a large degree Galois number field
to have root discriminant less than $\Omega$.  In the four cases of Table~\ref{fourtable}, 
the Galois root discriminants are quite close to $\Omega$.

\subsubsection{Comparison with fields from the Ramanujan newform.}
An alternative approach to keeping root discriminants small is simply to insist that levels $N$ 
be just $1$.   The smallest weight newform with $N=1$ is the famous Ramanujan form
$\Delta_{12} \in S_{12}(1)$.   Its projective mod $\ell$ Galois respresentations are known to be
surjective onto $PGL_2(\F_\ell)$ except for $\ell \in \{2,3,5,7,23,691\}$ \cite{SwD73}.   
For surjective representations with $\ell < 3500$, Elkies and Atkin showed there is no companion form of type 1T \cite[\S17]{gr}.  The first two 
$\ell$ for which there is a companion form of type 2T are known to be 
$2411$ and $7758337633$ \cite{LR10}.  When $\ell>7$ and no companion forms are present,  the slope content at $\ell$ is
$[(\ell+10)/(\ell-1)]_{(\ell-1)/w}$; here $w=10$ if $\ell \equiv 1 \; (11)$ and otherwise $w=1$.  
The resulting quantities, assuming $w=1$ in the case of $\Delta$, are
\begin{align*}
\delta & = \ell^{(\ell+10)/(\ell+1)}, & 
\Delta & = \ell^{(\ell^2+10 \ell - 12)/(\ell^2 - \ell)}.
\end{align*}
In the four examples of Table~\ref{fourtable}, the field from $\Delta_{12}$ is substantially
more ramified than the tame field.  

\subsubsection{Comparison with other fields}
Schaeffer found a weight $1$ modular form with level $N=3 \cdot 227$ and quadratic character $\chi_{-227}$ 
living in characteristic $11$ \cite[Table~4]{Sch15}.  
From this modular form, one knows that there is a $PGL_2(\F_{11})$ field
with root discriminant $\delta = 3^{10/12} 227^{5/12} \approx 23.94$ and Galois root discriminant
$\Delta = 3^{10/11} 227^{1/2} \approx 40.90$.  These quantities are much smaller than the corresponding
quantities on the $\ell = 11$ line of Table~\ref{fourtable}.    At present, one does not have a polynomial for this remarkable field.  

More explicitly, from an elliptic curve with conductor $128$ and $j$-invariant also $128$, one gets a 
degree fourteen polynomial with Galois group $PGL_2(\F_{13})$ and slope content 
$[8/3,8/3]_3$ at $2$ and $[13/12]_{12}$ at $13$.  This gives $\delta = 2^{13/7} 3^{13/14} \approx 39.21$ which
undercuts $43.00$ from Table~\ref{fourtable}.  However the Galois root discriminant $\Delta = 
2^{13/6} 13^{167/156} \approx 69.94$ 
is well above $47.82$.

 \subsection{The tabulation problem and group-drop}
 \label{tabulation}
 A standard problem in the theory of number fields goes as follows:  {\em Let $G$ be a transitive permutation
 group of degree $n$;  let $B$ be a positive real number;  determine the complete list 
 of degree $n$ number fields $K$ with associated Galois group $G$ and root discriminant $\leq B$.}   
 Often one thinks in terms of the ordered list of root discriminants, 
 $\delta_1(G) \leq \delta_2(G) \leq \delta_3(G) \leq \cdots$, with 
 particular interest in finding $\delta_1(G)$ for as many permutation groups 
 as possible.   
 
 The database
 \cite{jr-global-database} contains solutions of this problem for many small $G$ and large $B$.  
 For solvable groups $G$, class field theory lets one obtain non-empty lists for 
 quite large $G$ in quite large degree $n$.   However 
 for almost simple non-solvable groups realized in their lowest degree, for example $S_n$ itself or  
 $PGL_2(\F_\ell) \subseteq S_{\ell+1}$ for $\ell \geq 7$, 
 the standard purely number-theoretic approach rapidly decays from easy to impossible 
 as $n$ increases from $5$ to $10$.   The numbers  $\delta \approx 44.07$ and $49.50$ 
 from Table~\ref{four}
 are currently candidates for $\delta_1(G)$ with $G = PGL_2(\F_\ell) \subset S_{\ell+1}$, with  $\ell = 19$ and $29$
 respectively.   Similarly, the numbers $\Delta \approx 47.82$, $46.43$, and $50.62$ are candidates
 for $G = PGL_2(\F_\ell) \subset S_{\ell^3-\ell}$ and $\ell = 13$, $19$, and $29$ respectively.   
  
 For larger permutation groups $G \subseteq S_n$, as studied especially in \cite{kluners-malle}, 
 one often starts 
 with a parameterized family of fields $K_t$ having Galois group in
 $G$ for all $t$ and equal to $G$ for most $t$.    One searches
 among these fields for the $K_t$ with particularly small root discriminant.
 Often one encounters the phenomenon of {\em group-drop}: the 
 $K_t$ with the smallest root discriminant all have Galois group
 strictly smaller than $G$.  When this phenomenon occurs
 with great strength, it is some heuristic indication that the smallest
 root discriminant of a $G$-field found is not too far 
 above the actual minimum $\delta_1(G)$ sought. 
 
 The phenomenon of group-drop has indeed occurred with great strength behind the scenes 
 in \S\ref{mathieu}, \S\ref{obtaining}, and \S\ref{Polynomials1}.   We add some perspective now by 
 discussing some numerics of the group-drops observed.  
 
 \subsubsection{Drops from Mathieu groups to $PGL_2(\F_{11})$ in \S\ref{mathieu}.} 
 It seems that all the $M_{22}.2$ fields and $M_{12}.2$ fields obtain by generic 
 specializations of the covers mentioned in \S\ref{mathieu} are much more heavily ramified than the exceptional
 $PGL_2(\F_{11})$ specialization we are pursuing.  While $K_1$ has Galois root discriminant
 $\Delta_1 \approx 52.75$, the smallest GRD for an $M_{22}.2$ field that we 
 have found is $\Delta^{\rm gen}_1 = 2^{139/48} 3^{3/4} 11^{2/3} \approx 83.91$, coming from 
 $y = -5^2 139^2/2^7$ in the cover of \cite{ma}.  Similarly, while $K_2$ has $\Delta_2 \approx 58.55$, 
 the smallest GRD we have seen for an $M_{12}.2$ specialization is $\Delta_1^{\rm gen} = 2^{2/3} 3^{31/18} 11^{11/12} \approx 94.84$,
 coming from $t = -5^3/2^2$ \cite[Table~5.2]{ro-mathieu}.  
 
 \subsubsection{Drops from $PGL_2(\F_{\ell})$ to solvable groups in \S\ref{obtaining}}
 \label{twoquartics}
   In the process of
 searching for our fifty-three triples $(K,g,h)$, we encountered other triples $(K,g,h)$ 
 which satisfy all the required conditions, except that the image of the common projective 
 representation is not all of $PGL_2(\F_\ell)$.  
  Two such triples are reported via the $(8)$'s and $(4)$'s appearing
 in Figure~\ref{Type2}.   The Galois root discriminants, 
 calculated by the general formulas presented in \S\ref{ramificationmod},
  are $2^{7/6} 11^{3/4} \approx 13.56$ and 
  $2^{2/3} 19^{3/4} \approx 14.45$.  These numbers are so small that they contradict the unconditional
  lower bounds for fields of degree $|PGL_2(\F_{11})|$ and $|PGL_2(\F_{19})|$ 
  respectively.  In fact, the Galois group is $S_4$ in each case, defining
  polynomials being respectively $x^4 - 2 x^3 - 4 x^2 - 6 x - 2$ and $x^4 - x^3 - 2x^2 - 6x - 2$.
 
  

 \subsubsection{Drops from $PGL_2(\F_{11})$ to solvable groups in \S\ref{Polynomials1}.}
 As explained in \S\ref{Polynomials1}, the field $K_2$ from Section~\ref{two} arises 
 also from specialization of the Atkin modular polynomial for $\ell = 11$.  
A computer search shows at least $394$ values
 of $j$ which keep ramification within $\{2,3,11\}$.   The
 smallest seven GRDs are $20.70$, $24.48$, $24.77$, $25.48$, $29.34$, $32.45$,
 and $49.50$.  All of them come from degenerate specializations
 with solvable Galois group.  The remaining $387$ points all give
 Galois group $PGL_2(\F_{11})$, with the smallest GRD being $\Delta_2 \approx 58.55$ 
 from $K_2$.

 \subsection{Modular approaches to the tabulation problem}
 \label{modularapproaches}
 Let $\lambda$ be a power of a prime $\ell$.  Our concluding point is that 
 modular methods can be brought to bear on the tabulation problem for $G$ any subquotient
 of $GL_2(\F_\lambda)$, in any transitive permutation representation.    Outside of 
 $\ell=2$, modular methods do not see totally real fields.  However there
 are analytic lower bounds on the minimal root discrimininants of totally real fields in 
 degree $n$.  As $n \rightarrow \infty$, these 
 lower bounds tend unconditionally to $4 \pi  e^{ \gamma + 1} \approx 60.84$ \cite[(2.5)]{od-survey}. 
  When
 cutoffs are kept small enough, totally real fields are not present.  In fact, 
 under the generalized Riemann hypothesis, the above limit increases
 to $8 \pi e^{\gamma + \pi/2} \approx 215.33$ \cite[(2.6)]{od-survey}; so one does not
 expect to see totally real fields towards the beginning of tables at all. 
 
 In approaching this problem, it is natural to break into 
 three cases: $\ell$ is wildly ramified, $\ell$ is tamely ramified, 
 and $\ell$ is unramified.   In all cases, one needs to search among eigenforms
 without any rationality condition imposed, so that  general characters and thus odd weights $k$ are 
 considered as well.   As an example that stays mostly in the context of this paper, 
 the unique newform in $S^{\rm new}_6(8)$ has an irrational
  companion form in $S^{\rm new}_{20}(8)$ for $\ell=23$; this yields a field with Galois group 
  $PGL_2(\F_{23})$ and the small Galois root discriminant $\Delta = 2^{7/6} 23^{23/24} \approx 45.30$.
 
 For a given cutoff $B$, it should be easiest to obtain
 complete lists in the wild-at-$\ell$ case, as the levels to be considered would be very small.   
 Next easiest would be the tame-at-$\ell$ case, as modeled by our computations throughout
 this paper, including the previous paragraph;   
  since ramification at $\ell$ is lighter, the levels to be considered
  would no longer all be so small.   By far the hardest, with our current theoretical
  knowledge, would be the unramified-at-$\ell$ case.  This case requires computations
  either with weight $1$ or weight $\ell$ forms, both difficult for different reasons;
  also since there is no ramification at $\ell$, large levels would 
  have to be inspected.

\bibliographystyle{amsalpha}
\bibliography{jr}

 \end{document}